\documentclass[preprint,12pt,3p]{elsarticle}

\usepackage{graphicx}
\usepackage{graphics}
\usepackage{epstopdf} 
\usepackage{epsfig} 
\usepackage{mathptmx} 
\usepackage{times} 
\usepackage{amsmath} 
\usepackage{amssymb}  
\usepackage{amsmath}
\usepackage{enumerate}
\usepackage{mathrsfs}
\usepackage{multirow}
\usepackage{subfigure}
\usepackage{latexsym}
\usepackage{bm}
\usepackage{color}
\usepackage{url}
\usepackage[utf8]{inputenc}
\usepackage[english]{babel}
\usepackage{amssymb}
\usepackage[english]{babel}
\usepackage[utf8]{inputenc}
\usepackage{algorithm}
\usepackage[noend]{algpseudocode}
\usepackage{algorithmicx}
\usepackage{soul}

\newtheorem{lemma}{Lemma}[section]
\newtheorem{definition}{Definition}[section]
\newtheorem{assumption}{Assumption}[section]
\newtheorem{theorem}{Theorem}[section]

\newtheorem{remark}{Remark}[section]

\begin{document}
\begin{frontmatter}

\title{Hierarchical Decentralized Reference Governor using Dynamic Constraint Tightening for Constrained Cascade Systems} 

\author[label1]{Shahram Aghaei\corref{cor1}}
\address[label1]{ Department of Electrical Engineering, Yazd University, Yazd, Iran}
\cortext[cor1]{Corresponding author}
\ead{aghaei@yazd.ac.ir}

\author[label2]{Abolghasem Daeichian}
\address[label2]{Department of Electrical Engineering, Faculty of Engineering, Arak University, Arak, 38156-8-8349 Iran}

\author[label3]{Vicenç Puig}
\address[label3]{Department of Automatic Control, Universitat Politècnica de Catalunya (UPC), Barcelona, Spain}
\footnote{This is the peer reviewed version of the following article: Aghaei, Shahram, Abolghasem Daeichian, and Vicenç Puig. "Hierarchical decentralized reference governor using dynamic constraint tightening for constrained cascade systems." Journal of the Franklin Institute 357.17 (2020): 12495-12517., which has been published in ?final form at
DOI: https://doi.org/10.1016/j.jfranklin.2020.09.040. This article may be used for non-commercial purposes in accordance with ScienceDirect Terms and Conditions for Use of Self-Archived Versions.}

\begin{abstract}
This paper proposes a hierarchical decentralized reference governor for constrained cascade systems. The reference governor (RG) approach is reformulated in terms of receding horizon strategy such that a locally receding horizon optimization is obtained for each subsystem with a pre-established prediction horizon. The algorithm guarantees that not only the nominal overall closed-loop system without any constraint is recoverable but also the state and control constraints are satisfied in transient conditions. Also, considering unfeasible reference signals, the output of any subsystem goes locally to the nearest feasible value.
The proposed dynamic constraint tightening strategy uses a receding horizon to reduce the conservatism of conventional robust RGs. Moreover, a decentralized implementation of the algorithms used to compute tightened constraints and output admissible sets is introduced that allow to deal with large scale systems. 
Furthermore, a set of dynamic constraints are presented to preserve recursive feasibility of distributed optimization problem. Feasibility, stability, convergence, and robust constraint satisfaction of the proposed algorithm are also demonstrated. The proposed approach is verified by simulating a system composed of three cascade jacketed continuous stirred tank reactors.
\bigskip


\end{abstract}
\end{frontmatter}

\section{Introduction}
Reference Governor (RG) control is one of the most efficient approaches to deal with systems subject to state and control constraints. The idea behind RG is quite simple: a regulator is synthesized without considering constraints, then a feasible reference signal is computed that guarantees the constraints are always satisfied. Many algorithms have been proposed to design the RG for linear and nonlinear systems \cite{gilbert1995discrete,bemporad1998reference,kalabic2011reference,kolmanovsky2014reference,di2015reference,zare2020}. Comprehensive reviews on reference governor strategies have been presented in \cite{kolmanovsky2014reference,garone2017reference}.
Receding Horizon Control (RHC) (also known as Model Predictive Control (MPC)) is another control technique which has gained a wide popularity in industry due to its flexibility in definition of control objectives and explicitly considering the operational constraints of the process variables. For an overview of RHC, its industrial applications, and main theoretical results, the reader is referred to \cite{muller2017economic,camacho2013model,mayne2000constrained,christofides2011networked,limon2010robust,betti2012anMPC,mendes2017practical,wang2017non,Ghenaati2019}.
Integrating RHC and RG is investigated by researchers in order to avoid constraints violation in transient due to set-point change, enlarging the domain of attraction, and replacing the inadmissible steady state reference inputs with the nearest admissible manipulated values. For instance, advanced RHC simply cascades RG and RHC \cite{de2011reference,di2015coordinating,falugi2015model}, while, some others merge RHC and RG in a single control scheme \cite{shahram2013mpc,Mayne2016generalized} which guarantee both the performance recovering without considering constraints as in RG while providing the flexibility properties as in RHC.

Recently, decentralized and distributed RG strategies have been developed for large-scale systems subject to local and global constraints  (\cite{tedesco2012distributed,kalabic2013decentralized,casavola2014distributed,casavola2016parallel,garone2017reference,casavola2018distributed} and references therein). For instance, \cite{garone2009distributed} presents a static decentralized reference management problem which leads to two approaches: Sequential approach \cite{casavola2011distributed} and parallel approach \cite{casavola2011bdistributed}. The parallel approach works unsatisfactorily when the command is close to its boundary. In sequential distributed strategies, only one agent (or subsystem) per decision time is allowed to modify its command, while, in the parallel approach based distributed strategies, all the agents (or subsystems) are allowed to modify their own reference signals simultaneously, assuming that the other agents will make the worst-case choices.  Some hybrid approaches switching between parallel and sequential modes has also been suggested \cite{tedesco2012distributed}. Regarding the mentioned interesting approaches, the following points are noticeable: 
\begin{itemize}
\item To the best of our knowledge, the prediction horizon in RG approaches is usually  considered to be equal to one.  In  \cite{gilbert2011constrained}, an Extended Command Governor (ECG) algorithm is proposed showing that a prediction horizon more than two does not affect the performance. Additionally, not only dimension of the MOAS is depending on the prediction horizon but also the new MOAS must be recalculated for any considered new value of prediction horizon. 
\item In the presence of not measurable but bounded disturbances, CG strategies attempt to  govern the reference input with regard to the nominal model subject to the minimal Robust Positively Invariant (mRPI) sets which limit the domain of manipulated reference at the maximum of conservatism \cite{garone2017reference,casavola2014distributed,garone2011sensorless}. The use of tubes with dynamic constraint tightening  has proposed for first time in \cite{mayne2005robust} for robust model predictive control of constrained, linear, discrete-time systems in the presence of bounded disturbances which reduces the conservatism of conventional mRPI based MPC approaches. But, based on the authors knowledge, such a solution has not yet presented for RG strategies.
\item Distributed RGs, that have already been presented in \cite{casavola2014distributed,casavola2016parallel,kalabic2013decentralized,tedesco2012distributed}, updates the reference input only for one subsystem at each time step while the increment is limited on a specific static set. The reference inputs of other subsystems are frozen to their recent values. In other words, using the mentioned distributed RG strategies, the variation of manipulated reference input enforced to fulfill a static constraint set while determining that such a set is problematic. This fact increases the conservatism of the method also for nominal disturbance-free linear systems.
 \item  The aggregated reference vector is designed to belong to overall static Maximal Output Admissible Set (MOAS) while the centralized computation of such a set may not be achievable due to large dimension of system. 
\end{itemize}.

This paper introduces a hierarchical decentralized RG for constrained systems.  This family of systems appear, as e.g., in the study of irrigation, drainage and potable water networks \cite{Cantoni2007}. Cascade systems presents a hierarchical lower block triangular (LBT) structure as described in Siljak  \cite{siljak2011decentralized} and Lunze \cite{Lunze1992}. As discussed in these references, the LBT structure allows a sequential optimization from the top to the bottom of the hierarchy. The approach proposed in this paper benefits from this fact. \\

Taking into account the state of art regarding RG approaches presented above, the main contributions of this paper are:
\begin{itemize}
\item A hierarchical Decentralized Reference Governor (DRG) algorithm is proposed exploiting the LBT structure of the system that involves solving a sequence of Receding Horizon Optimization Problems (RHOPs) using the MPC-based RG approach, proposed in \cite{shahram2013mpc}, that calculates the references for the local RG corresponding to each subsystem. The solution of the RHOPs in a hierarchical order allows to compute the optimal feasible reference values. 
\item Since the DRG approach is based on the MPC-based RG approach proposed in \cite{shahram2013mpc}, the reference prediction horizon for the local RGs could be selected more than one. In  \cite{shahram2013mpc}, it has been shown that using this approach, the dynamic variation of manipulated reference along  the prediction horizon leads to faster closed-loop transient response compared to the standard RG based on one-step reference horizon.
 \item A new dynamic constraint tightening approach is proposed for the   MPC-based  RG  approach used in the DRG scheme to provide maximum possible variation for feasible reference signal of subsystems at any decision time. Compared to the static tightening approach used in \cite{shahram2013mpc}, the reference signal is updated by the maximum possible value of the feasible increment at each time step for all subsystems allowing to achieve better performance.
 \item A method for computing all necessary invariant sets, namely RPI sets and MOASs, locally in each subsystem profiting the LBT structure of the system is provided. This method avoids computing them a centralized manner that would not be feasible in case of large scale system because of the high dimensional invariant sets.
 \item The stability, feasibility, and convergence properties of the proposed DRG scheme are established by resorting to the analysis tools developed for robust RHC, see e.g. \cite{mayne2005robust} and exploiting the LBT structure of the system \cite{siljak2011decentralized}. 
 \end{itemize} 
   
The effectiveness of the proposed approach is validated in simulation using a three jacketed continuous stirred tank reactors (CSTR) that are connected in cascade.

The paper is organized as follows. In Section \ref{ProblemStatement}, the problem statement is presented introduced both the open-lood and closed model of the subsystems. Section \ref{RHRG Algorithm} presents the proposed dynamic contraint tightening  based-DRG. The problem is reformulated in terms of auxiliary representations of subsystems and some design models are presented. The local RHOPs and the output admissible sets for the corresponding nominal closed-loop subsystems are introduced. Finally, the local RHOPs are stated using the dynamic constraint tightening DRG algorithm. The stability, feasibility, and convergence properties  are also investigated. The proposed  DRG algorithm is applied to the cascade connection of three CTSR linearized models in Section \ref{SimulationExample}. Finally, the paper is concluded in Section \ref{conclusion}.

\textbf{Notations:} $F^T$ indicates the transpose of matrix $F$. $\oplus$ an $\ominus$ denote the Minkowski sum and Pontryagin difference, respectively (see \cite{mayne2005robust}) and $\bigoplus_{i=1}^{M}\mathbb{O}^{i}=\mathbb{O}^{1}\oplus...\oplus\mathbb{O}^{M}$. $\mathfrak{B}^{dim}_{\epsilon}(0)$ defines a ball of radius $\epsilon$ centered at the origin in the space $\mathbb{R}^{dim}$.

\section{Problem Statement} \label{ProblemStatement}
\subsection{Problem Definition}\label{ss211}

Let consider a large scale system  modeled using the following discrete-time state space model

\begin{eqnarray}\label{eq1}
x(k+1) &=&A x(k)+B u(k)+E w(k)  \nonumber\\
y(k)   &=&C z(k)
\end{eqnarray}
where $x(k)\in\mathbb{R}^{n_x}$, $u(k)\in\mathbb{R}^{n_u}$,  and $y(k)\in\mathbb{R}^{n_y}$, are state, input,  and controlled output vectors at time step $k$, respectively.  The disturbances $w(k)\in\mathbb{W}\subset\mathbb{R}^{n_w}$ are assumed to be unknown but bounded in the set $\mathbb{W}$. The states and inputs $x(k)\in\mathbb{X}\subset\mathbb{R}^{n_x}$  and $u(k)\in\mathbb{U}\subset\mathbb{R}^{n_u}$ are physically constrained.

As announced in the introduction, the system \eqref{eq1} is assumed to be composed by a cascade of subsystems. Because of this particular structure, the system may be transformed into hierarchical LBT decomposition \cite{siljak2011decentralized}, i.e.,  the system matrices $A$,  $B$ and $E$ a expressed as lower block triangular matrices. 

\begin{figure}[!h]
	\begin{center}
	\includegraphics[width=5cm]{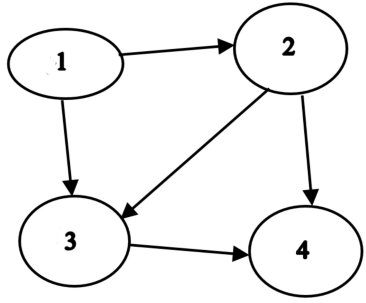}
	\caption{Cascade system}
	\label{fig:graph}
	\end{center}
\end{figure}

Before introducing local subsystem equations,  the concepts of inlet and outlet neighbor sets are recalled using graph theory \cite{godsil2013algebraic}:
\begin{definition}[Inlet neighbor set]	\label{D1}
	A set of all subsystems defined in \eqref{Subsystem} which are  affecting, directly, the $i^{th}$ subsystem is called the inlet neighbor set of the $i^{th}$ subsystem, i.e., $$\mathcal{N}_{IN}^{i}=\lbrace j<i\vert \Phi_{ij}\neq\textbf{0} \rbrace.$$
\end{definition}
\begin{definition}[Outlet neighbor set] \label{D2}
	A set of all subsystems defined in \eqref{Subsystem} which are affected, directly, by the $i^{th}$ subsystem is defined the outlet neighbor set of the $i$-th subsystem, i.e., $$\mathcal{N}_{OUT}^{i}=\lbrace j>i\vert \Phi_{ji}\neq\textbf{0} \rbrace.$$
\end{definition} 
In these definitions, the value matrix $\Phi_{ij}$ indicates an influencing matrix describing the direct influence of node $i$ on node $j$ which will be equal to matrix $\textbf{0}$ in cases of null interaction from node $i$ to node $j$. 
As an illustrative example of these definitions consider the system depicted in Figure \ref{fig:graph}. The cascade inlet and outlet neighbor sets are: $\mathcal{N}_{IN}^{1}=\lbrace\rbrace$, $\mathcal{N}_{IN}^{2}=\lbrace 1 \rbrace$, $\mathcal{N}_{IN}^{3}=\lbrace 1,2 \rbrace$, $\mathcal{N}_{IN}^{4}=\lbrace 2,3 \rbrace, \mathcal{N}_{OUT}^{1}=\lbrace2,3\rbrace$, $\mathcal{N}_{OUT}^{2}=\lbrace 3,4 \rbrace$, $\mathcal{N}_{OUT}^{3}=\lbrace 4 \rbrace$, and $\mathcal{N}_{OUT}^{4}=\lbrace\rbrace$.

Definition \ref{D1} is used to describe, in a more precise way, the set of $M$ discrete-time linear time-invariant (DLTI) subsystems of system \eqref{eq1} as follows: 
\begin{eqnarray}\label{Subsystem}
x^{i}(k+1) &=& A_{ii} x^i(k)+\sum_{j\in\mathcal{N}_{IN}^{i}}(A_{ij}x^{j}(k))+B_{ii} u^i(k)+E_{ii} w^{i}(k)  \nonumber\\
y^i(k)   &=& C_{ii} x^i(k)  \nonumber\\
c^i(k)   &=& \left[x^i(k)^T,u^i(k)^T\right]^T, i=1,..,M,
\end{eqnarray}
where $x^i(k)\in\mathbb{R}^{n_i}$, $u^i(k)\in\mathbb{R}^{p_i}$, $w^i(k)\in\mathbb{W}^{i}\subset\mathbb{R}^{n_{w_i}}$, $y^i(k)\in\mathbb{R}^{p_i}$, and $c^{i}(k)\in\mathbb{R}^{n_i+p_i}$ are state, input,  disturbance, output, and constraint vectors at time step $k$, respectively. Also, $\sum_{i=1}^{M}n_{x_{i}}=n_x$, $\sum_{i=1}^{M} n_{u_{i}}=n_u$, and $\sum_{i=1}^{M} n_{w_{i}}=n_w$. Constraint sets for the overall system variables in \eqref{Subsystem} satisfy $\prod_{i=1}^{M}\mathbb{X}^{i}=\mathbb{X}$, $\prod_{i=1}^{M}\mathbb{U}^{i}=\mathbb{U}$, and $\prod_{i=1}^{M}\mathbb{W}^{i}=\mathbb{W}$. The sets $\mathbb{X}^{i}$, $\mathbb{U}^{i}$, and $\mathbb{W}^{i}$ are supposed to be closed, compact and contain the origin for all $i=1,\cdots,M$.

Definition \ref{D2} will be used in Section \ref{RHRG Algorithm} for formalizing the proposed algorithms.

\subsection{Local Controllers and Assumptions}\label{ss212}
Assume that each subsystem is stabilized by a local controller that guarantees good tracking performance. So, $i^{th}$-subsystem has local reference input $\check{g}^i(k)\in\mathbb{R}^{p_i}$ and local controlled output $y^i(k)\in\mathbb{R}^{p_i}$ that evolves according to

\begin{eqnarray}\label{InterconnectedSubsystem}
z^i(k+1) &=& \Phi_{ii}  z^i(k)+\sum_{j\in\mathcal{N}_{IN}^{i}}\Phi_{ij} z^j(k)+\Gamma_{ii}\check{g}^i(k)+\Omega_{ii} w^i(k) \nonumber\\
y^i(k)     &=& \Upsilon_{ii} z^i(k) \nonumber\\
c^i(k)     &=& H_{ii} z^i(k)
\end{eqnarray}

\noindent where $z^i(k)\in\mathbb{R}^{n_{z_i}}$ for $i=1,\cdots,M$ represent locally closed-loop state vectors at time step $k$. Also, regarding the open-loop subsystems interaction matrix $A_{ij}$ for any $j\in\mathcal{N}_{IN}^{i}$, $\Phi_{ij}$ describes the direct influence of any locally closed-loop subsystem $j$ affecting the state of locally closed-loop subsystem $i$. Note that $\Upsilon_{ii} z^i(k)=C_{ii} x^i(k)$ and $H_{ii} z^i(k)=\left[x^i(k)^T,u^i(k)^T\right]^T$. So, the corresponding constraints that have to be fulfilled are $x^{i}\in\mathbb{X}^{i}$ and $u^{i}\in\mathbb{U}^{i}$, i.e. $c^i(k)\in\mathbb{X}^{i}\times\mathbb{U}^{i}$.

\begin{assumption} \label{A1}
	For each subsystem $i$:
	\begin{enumerate}
		\item The state vector is measurable;
		\item $\Phi_{ii}$ is a Schur matrix;
		\item $\Upsilon_{ii}(I-\Phi_{ii})^{-1}\Gamma_{ii}=I_{n_{y_{i}}}$
		\item $\Upsilon_{ii}(I-\Phi_{ii})^{-1}\Phi_{ij}=0_{n_{y_i}\times n_{x_j}}, \quad \forall j\in\mathcal{N}_{IN}^{i}$.
	\end{enumerate}
\end{assumption}
This assumption guarantees that at least \emph{in nominal conditions}, i.e. when the state and control constraints are not active and the disturbances $w^i$ are null, each subsystem \eqref{InterconnectedSubsystem} ensures offset-free tracking performance regarding its local constant reference input. The conditions 3 and 4 in Assumption \ref{A1} are standard specifications in feedback control design (see e.g. \cite{Franklin1997}). Both could be hold simultaneously according to the following lemma.

\begin{lemma}\label{L1}
Consider that each subsystem  \eqref{Subsystem} has been locally augmented with tracking error integrators. Then, by designing an inner loop state feedback regulator, conditions 3 and 4 in Assumption \ref{A1} hold, simultaneously, in the steady state response to the constant reference inputs.
\end{lemma}

\noindent\textbf{Proof:}\\

To achieve a good tracking performance in case of constant reference inputs, according to the internal model principle, dimensional compatible error integrator dynamics could be applied as the controller and then the augmented system/model could be stabilized by properly designed state feedback. Proceeding in such a way for each subsystem $i$, the local model of any subsystem augmented by error integrating controller can be written as follows
 
\begin{align} \label{eq:augmented}
  & {{z}^{i}}(k+1)={{{\bar{\Phi }}}_{ii}}{{z}^{i}}(k)+\sum\limits_{j\in N_{IN}^{i}}{{{\Phi }_{ij}}{{z}^{j}}(k)}+{{\Gamma }_{ii}}{{u}^{i}}(k)+{{\Omega }_{ii}}w(k) \\ 
 & {{y}^{i}}(k)={{\Upsilon }_{ii}}{{z}^{i}}(k) \nonumber 
\end{align}

\noindent with

\begin{align}
{{z}^{i}}(k)&=\left[ \begin{matrix}
   {{x}^{i}}(k)  \\
   x_{a}^{i}(k)  \\
\end{matrix} \right],\quad{{\bar{\Phi }}_{ii}}=\left[ \begin{matrix}
   {{A}_{ii}} & {{0}_{{{n}_{x_i}}\times {{n}_{y_i}}}}  \\
   -{{C}_{ii}} & {{I}_{{{n}_{y_i}}}}  \\
\end{matrix} \right],\quad {{\Phi }_{ij}}=\left[ \begin{matrix}
   {{A}_{ij}}  \\
   {{0}_{{{n}_{y_i}}\times {{n}_{x_j}}}}  \\
\end{matrix} \right], \quad \\{{\Gamma }_{ii}}&=\left[ \begin{matrix}
   {{B}_{ii}}  \\
   {{0}_{{{n}_{y_i}}\times {{n}_{y_i}}}}  \\
\end{matrix} \right],\quad{{\Omega }_{ii}}=\left[ \begin{matrix}
   {{0}_{{{n}_{x_i}}\times {{n}_{{{w}_{i}}}}}}  \\
   {{E}_{ii}}  \\
\end{matrix} \right],\quad {{\Upsilon }_{ii}}=\left[ \begin{matrix}
   {{C}_{ii}} & {{0}_{{{n}_{y_i}}\times {{n}_{y_i}}}}  \\
\end{matrix} \right],\nonumber 
\end{align}
\noindent where $x_{a}^{i}$ is the state describing the integrating controller dynamics which dimension is equal to the number of local controlled outputs ${n}_{y_i}$.
 Let the following state feedback be designed as the stabilizing controller for the augmented system:
\[{{u}^{i}}(k)=\left[ \begin{matrix}
   {{K}_{{{x}^{i}}}} & {{K}_{x_{a}^{i}}}  \\
\end{matrix} \right]{{z}^{i}}(k),\]
Then, by applying this  state feedback  controller to the augmented system \eqref{eq:augmented},  the closed-loop system dynamics are governed by
\begin{align}
  & {{z}^{i}}(k+1)={{\Phi }_{ii}}{{z}^{i}}(k)+\sum\limits_{j\in \mathcal{N}_{IN}^{i}}{{{\Phi }_{ij}}{{z}^{j}}(k)}+{{\Gamma }_{ii}}\check{g}(k)+{{\Omega }_{ii}}w(k) \\ 
 & {{y}^{i}}(k)={{\Upsilon }_{ii}}{{z}^{i}}(k)  
\end{align}
\noindent where
\begin{align}{{\bar{\Phi }}_{ii}}={{\bar{\Phi }}_{ii}}-{{\Gamma }_{ii}}K=\left[ \begin{matrix}
   {{A}_{ii}}-{{B}_{ii}}{{K}_{{{x}^{i}}}} & -{{B}_{ii}}{{K}_{x_{a}^{i}}}  \\
   -{{C}_{ii}} & {{I}_{{{n}_{y_{i}}}}}  \nonumber
\end{matrix} \right]
\end{align}

So, the system state and output in function of the reference input $\check{g}$ , in the steady state, is given by
\begin{align}
  & {{z}^{i}}(\infty )={{\Phi }_{ii}}{{z}^{i}}(\infty )+\sum\limits_{j\in \mathcal{N}_{IN}^{i}}{{{\Phi }_{ij}}{{z}^{j}}(\infty )}+{{\Gamma }_{ii}}\check{g}(\infty ) \\ 
 & {{y}^{i}}(\infty )={{\Upsilon }_{ii}}{{z}^{i}}(\infty )={{\Upsilon }_{ii}}{{(I-{{\Phi }_{ii}})}^{-1}}{{\Gamma }_{ii}}\check{g}(\infty )+{{\Upsilon }_{ii}}{{(I-{{\Phi }_{ii}})}^{-1}}\sum\limits_{j\in \mathcal{N}_{IN}^{i}}{{{\Phi }_{ij}}{{z}^{j}}(\infty )}  
\end{align}

Applying the matrix inversion lemmas and considering the structure of the closed-loop system matrices, it yields

\begin{align}
  & {{\Upsilon }_{ii}}{{(I-{{\Phi }_{ii}})}^{-1}}=\left[ \begin{matrix}
   {{C}_{ii}} & 0  \\
\end{matrix} \right]\left[ \begin{matrix}
   {{0}_{{{n}_{x_i}}}} & {{(I-({{A}_{ii}}-{{B}_{ii}}{{K}_{{{x}^{i}}}}))}^{-1}}(-{{B}_{ii}}{{K}_{x_{a}^{i}}}){{(-{{C}_{ii}}{{(I-({{A}_{ii}}-{{B}_{ii}}{{K}_{{{x}^{i}}}}))}^{-1}}(-{{B}_{ii}}{{K}_{x_{a}^{i}}}))}^{-1}}  \\
   * & *  \\
\end{matrix} \right] \\ 
 & \quad \quad \quad \quad \quad \ =\left[ \begin{matrix}
   {{0}_{{{n}_{y_{i}}}\times {{n}_{x_{i}}}}} & {{I}_{{{n}_{y_{i}}}}}  \nonumber
\end{matrix} \right]  
\end{align}

\noindent from where follows the satisfaction of conditions 3 and 4 in Assumption \ref{A1}

\begin{align}
&{{\Upsilon }_{ii}}{{(I-{{\Phi }_{ii}})}^{-1}}{{\Gamma }_{ii}}={{I}_{{{n}_{y_{i}}}}}\\
&{{\Upsilon }_{ii}}{{(I-{{\Phi }_{ii}})}^{-1}}{{\Phi }_{ij}}=0_{n_{y_i}\times{n_{x_j}}}, \quad \forall j\in \mathcal{N}_{IN}^{i}
\end{align}

\textcolor{blue}{\hfill$\square$}

\subsection{Decentralized Reference Governor Objective}\label{ss213}

This paper addresses the DRG design problem for the subsystems \eqref{InterconnectedSubsystem} which consists of locally determining a feasible reference signal $\check{g}^i(k)$ to be the best approximation of main reference input $y^i_r(k)$ such that prevents from any constraint violation in $i$-th subsystem and guarantees $c^{i}(t)\in\mathbb{X}^{i}\times\mathbb{U}^{i}$, $\forall t>k$, $\forall i=1,\cdots,M$. To this end, we propose that the feasible reference input to each subsystem $i$ be the sum of main reference input $y^i_r(k)$ and a correction factor $\alpha^i(k)$:
\begin{equation}\label{reference}
\check{g}^i(k)=y^i_r(k)+\alpha^i(k)
\end{equation}
The DRG should determine the minimum value of $\alpha^i(k)$ which makes the $\check{g}^i(k)$ to be the nearest feasible reference input to $y_r^i(k)$. In a particular situation, $\alpha^i(k)$ should be zero when $y_r^i(k)$ is feasible in $i$-subsystem.


\section{Proposed Approach} \label{RHRG Algorithm}
This section presents the proposed approach: a hierarchical Dynamic Constrained Tightening Decentralized Reference Governor (DCT-DRG) algorithm for the closed-loop subsystems \eqref{InterconnectedSubsystem} based on sequential control strategy to find the best $\alpha^i(k)$ in Eq. \eqref{reference}. The method has been derived as a decentralized supervisory scheme to deal with state and control constraints by determining and applying, locally, feasible reference inputs. First, some dynamical models that are required in the design procedure and method for locally computing tightened constraint sets are presented. Then, an algorithm is introduced for the decentralized computation of all output admissible sets allowing to deal with large scale systems. Finally, the proposed DCT-DRG is presented where a novel \emph{dynamic constraint set} is introduced to preserve the recursive feasibility of locally optimization problems. At the end of this section, the recursive feasibility, stability and convergence properties of the proposed approach will be proven.

\subsection{Design models and tightened constraint sets}\label{ss31}
The following auxiliary design models are required for developing DCT-DRG (see Figure \ref{fig:scheme}):

\begin{figure}[htpb!]
	\begin{center}
		\includegraphics[width=10cm]{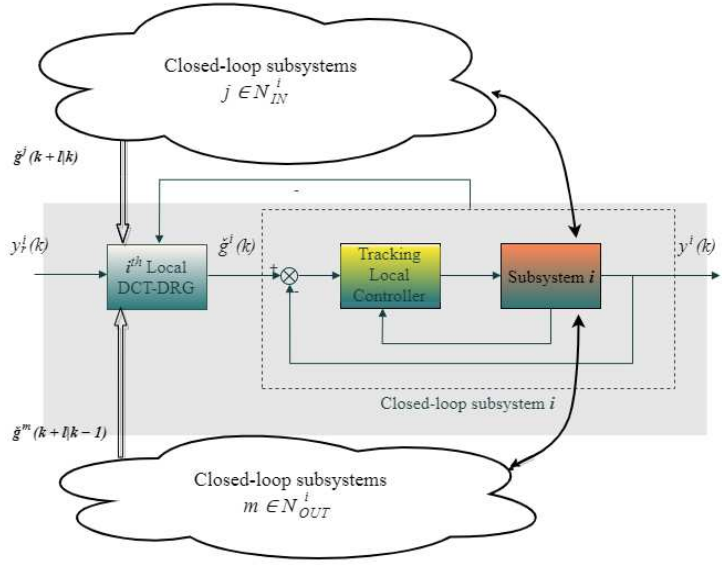}
		\caption{Schematic diagram of proposed DCT-DRG approach}
		\label{fig:scheme}
	\end{center}
\end{figure}

\begin{itemize}
\item[1)] \emph{Closed-loop Nominal Model} (CNM) which is used to predict the state trajectory in the prediction horizon. This model is obtained by neglecting the disturbance term in Eq. \eqref{InterconnectedSubsystem}:
\begin{eqnarray}\label{CNM}
z_c^i(k+1) &=& \Phi_{ii} z_c^i(k)+\sum_{j\in\mathcal{N}_{IN}^{i}}\Phi_{ij} z_c^j(k)
+\Gamma_{ii}(y^i_r(k)+\alpha^i(k)) \nonumber\\
y^i_c(k) &=& \Upsilon_{ii} z_c^i(k) \nonumber\\
c^i_c(k) &=& H_{ii} z_c^i(k).
\end{eqnarray}
 \noindent where $ z_c^i\in \mathbb{R}^{n_{z_i}}$ and $y^{i}_{c}\in \mathbb{R}^{n_{y_i}}$.

\item[2)] For calculating the so-called Robustly Positive Invariant (RPI) set for any subsystem \cite{rakovic2005invariant}, the \emph{Closed-loop Nominal Error Model} (CNEM) is introduced. This model is obtained by defining
the error signal $e^i= z^i-z_c^i$, and deriving its dynamical model by using \eqref{InterconnectedSubsystem} and \eqref{CNM}:
\begin{eqnarray}\label{CE}
	e^i(k+1) &=& \Phi_{ii}e^i(k)+ w_e^i(k).
\end{eqnarray}
where
\begin{eqnarray}\label{we}
	w_e^i(k) &=& \sum_{j\in\mathcal{N}_{IN}^{i}}\Phi_{ij}e^j(k)+\Omega_{ii} w^i(k)
\end{eqnarray}
The original constraint sets of Eq. \eqref{InterconnectedSubsystem} are reformulated in terms of CNM variables in two categories:
\\a) \textit{Transient tightened constraint set:} At any initializing time step $k$, by applying the measured state of subsystem $i$ as the initial state of CNM  \eqref{CNM}, i.e. $z_c^i(k)=z^i(k)$, one can write the constraint sets of Eq. \eqref{InterconnectedSubsystem} as:

\begin{eqnarray}\label{newconst}
c_c^i(k)=H_{ii} z^i_c(k) &\in& \mathbb{X}\mathbb{U}^{i}(k), \forall i=1,\cdots,M
\end{eqnarray}

where 

\begin{equation}\mathbb{X}\mathbb{U}^{i}(0)=\mathbb{X}^{i}\times \mathbb{U}^{i}
\end{equation}

with

$\mathbb{X}\mathbb{U}^{i}(k+1)=\mathbb{X}\mathbb{U}^{i}(k)\ominus H_{ii}\Phi_{ii}^k \mathbb{W}_{e}^{i}(k)$ and $\mathbb{W}_{e}^{i}(k)=\bigoplus_{j\in\mathcal{N}_{IN}^{i}}\Phi_{ij}\mathbb{W}_{e}^{j}(k)\oplus\Omega_{ii}\mathbb{W}^{i}$. $\mathbb{W}^{1}_{e}(0)=\Omega_{11}\mathbb{W}^{1}$ is given as the initialisation value.
\\b) \textit{Steady-state tightened constraint:} Since all $\Phi_{ii}$ for $i=1,...,M$ are Schur, then the corresponding transient tightened constraint sets converge to their steady state sets denoted as $\mathbb{X}\mathbb{U}^{i}_{\infty}=\mathbb{X}\mathbb{U}^{i}(\infty)$. The computation of such a sets are performed in terms of CNM variables by utilizing the RPI sets which are recursively calculated. In other words, $w_{e}^{1}(k)=\Omega_{11}w^{1}(k)\in\Omega_{11}\mathbb{W}^{1}=\mathbb{W}^{1}_{e}$ is given and a polytopic outer approximation $\mathbb{F}_{\infty}^{1}$ of the minimal RPI (mRPI) of Eq. \eqref{CE} can be computed with the method proposed in \cite{rakovic2005invariant}. So, by iterating for subsystems $i>1$, the RPIs $\mathbb{F}_{\infty}^{j}$, $j\in\mathcal{N}_{IN}^{i}$ could be computed through Eq. \eqref{we}. As a result, $w^{i}_{e}\in\mathbb{W}_{e}^{i}=\bigoplus_{j\in\mathcal{N}_{IN}^{i}}\Phi_{ij}\mathbb{F}_{\infty}^{j}\oplus\Omega_{ii}\mathbb{W}^{i}$. After calculating all the RPIs and considering $z_c^i= z^i-e^i$, it is easy to see that the original constraints can be reformulated in terms of the variables $z^i_c$ by defining the tightened constraint set
\begin{eqnarray} \label{TC}
\mathbb{X}{\mathbb{U}}^{i}_{\infty} &=& (\mathbb{X}^{i}\times \mathbb{U}^{i})\ominus
H_{ii}\mathbb{F}^{i}_{\infty}.
\end{eqnarray}
\begin{assumption}\label{ass:1}
${\mathbb{X}}{\mathbb{U}}^{i}(k)$ for any $i=1,...,M$, and any $k>0$ are closed polytopes which contain the origin in their interior.
\end{assumption}

\item[3)] \emph{Uncontrolled Closed-loop Nominal Model} (UCNM) which is obtained from \eqref{CNM} by shifting $\alpha^i$ one step backwards in time:
\begin{eqnarray}\label{UCNM}
z_u^i(k+1) &=& \Phi_{ii} z_u^i(k)+\sum_{j\in\mathcal{N}_{IN}^{i}}\Phi_{ij} z_u^j(k)+\Gamma_{ii}(y^i_r(k)+\alpha^i(k-1)) \nonumber\\
y_u^i(k)   &=& \Upsilon_{ii} z_u^i(k).
\end{eqnarray}
where $z_u^i\in\mathbb{R}^{n_{z_i}}$ and $y_u^i\in\mathbb{R}^{p_i}$. This is utilized as a reference model for the CNM.

\item[4)] \emph{Controlled Error Model} (CEM) which is obtained by defining $\varepsilon_d^i(k)= z_c^i(k)- z_u^i(k)$ and deriving its dynamical model as:
\begin{eqnarray} \label{EM}
\varepsilon^i_d(k+1) &=& \Phi_{ii} \varepsilon^i_d(k)+\Gamma_{ii}\delta\alpha^i(k)
\end{eqnarray}
where $\varepsilon_d^i\in\mathbb{R}^{n_{z_i}}$ and $\delta\alpha^i(k)=\alpha^i(k)-\alpha^i(k-1)$. This error should converge to zero at steady state.
\end{itemize}

\subsection{Decentralized Maximal Output Admissible set (MOAS)}\label{ss32}


The goal is to find a set where $z_c^i(k)$ must lie in order to guarantee that \eqref{newconst} are fulfilled for all $t\geq k$. Consider the CNM \eqref{CNM} with constant $\alpha^i$ and constant references $y^i_r$. We may write the model by defining $\tilde{y}^i_r=y^i_r+\alpha^i$ as:
\begin{subequations}
	\label{eq:exp}
	\begin{align}
	\begin{bmatrix}{ z}_c^i(k+1)\\\tilde{y}_r^i(k+1)
	\end{bmatrix}&=\mathcal{A}^{i}\begin{bmatrix}{ z}_c^i(k)\\\tilde{y}_r^i(k)
	\end{bmatrix}+\mathcal{B}^iw_{ z}^i(k)\label{eq:exp_state}\\
	c^{i}_{c}(k)&=\mathcal{C}^{i}\begin{bmatrix}{ z}_c^i(k)\\\tilde{y}_r^i(k)
	\end{bmatrix}
	\label{eq:exp_out}
	\end{align}
\end{subequations}
where
$$\mathcal{A}^i=\begin{bmatrix}\Phi_{ii} & \Gamma_{ii} \\
0 & I \end{bmatrix},\,\mathcal{B}^i=\begin{bmatrix}I\\0\end{bmatrix},\,
\mathcal{C}^{i}= \begin{bmatrix}
H_{ii} & 0\end{bmatrix}$$
and the interlacing term $w_z^i(k)=\sum_{j\in\mathcal{N}_{IN}^{i}}\Phi_{ij}z_c^j(k)\in\mathbb{W}^i_z$ is a bounded disturbance whose magnitude is known in advance (as specified later in the paper  in Algorithm 1). It is assumed that the pairs $(\mathcal{A}^{i},\mathcal{C}^{i})$ are observable for all $i=1,\cdots,M$.

The Maximal Output Admissible Set (MOAS) is a state invariant set such that the corresponding outputs satisfy prescribed constraints and is denoted by $\mathbb{O}^i_\varepsilon$. This set guarantees that if $(z_c^i(k),\tilde{y}_r^i(k))\in\mathbb{O}^i_\varepsilon$ then $(z_c^i(t),\tilde{y}_r^i(t))\in\mathbb{O}^i_\varepsilon$ and $H_{ii}z_c^i(t)\in\mathbb{X}\mathbb{U}^i_{\infty}(t), \forall t\geq k$. So, the set where $z_c^j(t)$ lies for all $t\geq k$ is ${\mathbb{O}}^{j}_{z}=\begin{bmatrix}I&0\end{bmatrix}{\mathbb{O}}^{j}_{\varepsilon}$.
Considering Assumption \ref{ass:1} and the observability of $(\mathcal{A}^{i},\mathcal{C}^{i})$, the following set is calculated for all subsystems by employing the algorithm proposed in \cite{KolGi98}
\begin{equation}\label{eq16}
\mathbb{O}^i_\varepsilon=\left\{(z_c^i,\tilde{y}_r^i) \mid \mathcal{C}^i(\mathcal{A}^i)^k(z_c^i,\tilde{y}_r^i)\in\mathbb{X}_\mathbb{U}^i(k), \forall k\geq 0, H_{ii}(I-\Phi_{ii})^{-1}\Gamma_{ii} \tilde{y}_r^i \in\mathbb{X}_{\mathbb{U},\varepsilon}^i\right\}
\end{equation}
where, by assuming $\mathbb{X}_\mathbb{U}^i(0)=\mathbb{X}\mathbb{U}^i_{\infty}$, $\mathbb{X}_\mathbb{U}^i(k)$ is computed by $\mathbb{X}_\mathbb{U}^i(k+1)=\mathbb{X}_\mathbb{U}^i(k)\ominus H_{ii}\Phi_{ii}^k \mathbb{W}^i_{z}$ for any $k>0$ and $\mathbb{X}_{\mathbb{U},\varepsilon}^i$ is the compact and convex set which contains the origin in its interior and arbitrarily close to $\mathbb{X}_\mathbb{U}^i(\infty)$, satisfying $\mathbb{X}_{\mathbb{U},\varepsilon}^i\subseteq\mathbb{X}_\mathbb{U}^i(\infty)\ominus \mathfrak{B}^{n_{d_i}}_{\varepsilon}(0)$; see \cite{KolGiACC95,KolGi98}. Also, the condition $H_{ii}(I-\Phi_{ii})^{-1}\Gamma_{ii}\tilde{y}_r^i\in\mathbb{X}_{\mathbb{U},\varepsilon}^i$ guarantees that $\lim_{k\rightarrow+\infty}\mathcal{C}^i(\mathcal{A}^i)^k(z_c^i(k),\tilde{y}_r^i)\in\mathbb{X}_{\mathbb{U},\varepsilon}^i$.
It is easy to verify that, if $(z_c^i(0),\tilde{y}_r^i(0))\in\mathbb{O}^i_\varepsilon$ then the constraints \eqref{newconst} are fulfilled for all $k\geq 0$ and $(z_c^i(t),\tilde{y}_r^i(t))\in\mathbb{O}^i_\varepsilon$, i.e., $\mathbb{O}^i_\varepsilon$ is RPI with respect to disturbances bounded in $\mathbb{W}^i_{z}$; see \cite{KolGi98}.
\begin{remark}	\label{Remark:R00}
 For computation of $\mathbb{O}^i_\varepsilon$, a practical alternative, more conservative with less computation effort, approach is based on using minimal RPI sets \cite{KolGiACC95,KolGi98} instead of maximal admissible ones by utilizing the RPI sets $\mathbb{X}_\mathbb{U}^i(\infty)$ in place of $\mathbb{X}_\mathbb{U}^i(k)$ in \eqref{eq16} for any $k$ and for any subsystem $i$.
\end{remark}

Now, the procedure for the decentralized computation of local MOAS can be formalized as presented Algorithm 1. This algorithm guarantees the fulfillment of \eqref{newconst} for all $t\geq k$.

\begin{algorithm}\label{AlMOAS}
	\caption{Decentralized Computation of local MOAS corresponding to any subsystem}
	\begin{algorithmic}
		\item[1)] \emph{Initialization}: compute $\mathbb{O}^1_{\varepsilon}$ by Eq. \eqref{eq16} considering $\mathbb{W}_z^1=\emptyset$ and $\mathbb{O}^1_z=\begin{bmatrix}I&0\end{bmatrix}{\mathbb{O}}^{1}_{\varepsilon}$
		\item[2)] \emph{Iteration}: For $i=2,\cdots,M$, repeat:\\
		\begin{itemize}
			\item[2-1)] Compute $\mathbb{W}_z^i= \bigoplus_{j\in\mathcal{N}_{IN}^{i}}\Phi_{ij}\mathbb{O}_z^j$.
			\item[2-2)] Compute $\mathbb{O}^i_{\varepsilon}$ by Eq. \eqref{eq16}.
			\item[2-3)] Compute $\mathbb{O}^i_z=\begin{bmatrix}I&0\end{bmatrix}{\mathbb{O}}^{i}_{\varepsilon}$.
		\end{itemize}
	\end{algorithmic}
\end{algorithm}

\begin{remark}	\label{Remark:R000}
As it is presented above, Algorithm 1 may fail in case ${\mathbb{O}}^{i}_{\varepsilon}=\emptyset$ for a given subsystem $i$. However, in case this happens, it is always possible to reduce $\mathbb{W}_{ z}^i$ (hence making problem feasible) by suitably enlarging the sets $\mathbb{W}_{ z}^j$, for all $j$ such that $\Phi_{ij}\neq 0$. This workaround has the effect of reducing the dimensions of sets ${\mathbb{O}}^{j}_{\varepsilon}$ for all the involved subsystems (i.e., the predecessors of $i$), but it does not hamper their invariance properties.
\end{remark}


\subsection{Decentralized Receding Horizon Optimization Problem} \label{ss33}
The main idea of the DCT-DRG algorithm relies on the solution of a set of $M$ standard Receding Horizon Optimization Problems (RHOP) formulated in terms of the variables $z_u^i$ and $\varepsilon_d^i$, or equivalently, $z_c^i=\varepsilon_d^i+z_u^i$. The RHOPs compute exogenous signal $\delta\alpha^i$ during the prediction horizon by minimizing a suitable performance index subject to a proper set of constraints, including \eqref{newconst}. However, at the end of prediction horizon $\delta\alpha^i=0$,  it is necessary to guarantee that the state and control constrains are still fulfilled. This can be achieved if states of the involved systems belong to suitably defined MOAS.

The proposed approach is now formulated for each subsystem $i$ in order to compute the augmented variable $\alpha^i$ of the proposed DCT-DRG. At any time instant $k$, the optimization problems must be solved sequentially from $i=1$ to $M$ considering that the solution of $j$-th RHOP influences the $i$-th RHOP where $j\in \mathcal{N}_{IN}^{i}$.

Let the length of prediction horizon to be $N\in\mathbb{N}$. Denote the predicted value of any variables $\chi$ at time step $k$ which is predicted at time step $t$ by $\hat{\chi}(k|t)$. Define
$\delta\hat\alpha^i(k:k+N|k) =(\delta\hat{\alpha}^i(k|k),\dots,\delta\hat{\alpha}^i(k+N|k))$. Then, optimal control sequence $\delta\hat\alpha^i(k:k+N-1|k)$ is given by solving the following optimization problem:

\begin{eqnarray}\label{eq18}
J_N^{i^{*}}\left(\varepsilon_d^i(k),y_r^i,\alpha^i(k-1)\right) &=&
\min_{\delta\hat\alpha^i(k:k+N-1|k)} J^i_N\left(\delta\hat\alpha^i(k:k+N-1|k);\varepsilon_d^i(k),y_r^i,\alpha^i(k-1)\right)
\end{eqnarray}
where
\begin{eqnarray}
J_N^i &=& \Vert\hat{\varepsilon}^i_{d}(k+N|k)\Vert_{P_{ii}}^{2}+\Vert\hat{\alpha}^i(k+N-1|k)\Vert_{P_{\alpha_{ii}}}^{2}+\sum_{l=0}^{N-1} \lbrace\Vert\hat{\varepsilon}^i_{d}(k+l|k)\Vert_{Q_{ii}}^{2}+\Vert\delta\hat{\alpha}^i(k+l|k)\Vert_{R_{\alpha_{ii}}}^{2}\rbrace\nonumber
\end{eqnarray}
subject to the constraints:
\begin{subequations}\label{eq:costrs}
\begin{align}
&H_{ii}\left(\hat{z_{u}}^{i}(k+l|k)+\hat{\varepsilon}^{i}_{d}(k+l|k)\right)\in\mathbb{X}{\mathbb{U}}^{i}(l), \qquad \forall l=1,\cdots,N-1, \label{eq:constr_stage}
\\
&H_{mm}\left(\hat{z}_c^m(k+l|k-1)+\hat{\sigma}^m_d(k+l|k)\right)\in\mathbb{X}{\mathbb{U}}^{m}(l+1), \qquad \forall m\in\mathcal{N}_{OUT}^i, l=1,\cdots,N-2, \label{eq:constr_aniticipative}
\\
&\sum_{j\in\mathcal{N}_{IN}^{m}}\Phi_{mj}(\hat{z}_c^j(k+N-1|k-1)+\hat{\sigma}^j_d(k+N-1|k))\in\mathbb{W}_{ z}^{m}, \qquad \forall m\in\mathcal{N}_{OUT}^{i}, \label{eq:terminal_sigma}
\\
&\begin{bmatrix}
\hat{z_{u}}^{i}(k+N|k)+\hat{\varepsilon}^{i}_{d}(k+N|k) \\
y_r^i+\hat{\alpha}^{i}(k+N-1|k)\end{bmatrix} \in\mathbb{O}^i_{\varepsilon}, \label{eq:constr_term}
\end{align}
\end{subequations}
where $P_{ii}$, $P_{\alpha_{ii}}$, $Q_{ii}$, and $R_{\alpha_{ii}}$ are weighting matrices. $\hat{z}_c^j(k+l|k)=\hat{z}_u^j(k+l|k)+\hat{\varepsilon}^j_d(k+l|k), \forall j\in\mathcal{N}_{IN}^i$ is the optimal prediction of the variable transmitted by subsystem $j$ to subsystems $i$. $\hat{\alpha}^i(k+l|k)$, $\hat{\varepsilon}^{i}_{d}(k+l|k)$, $\hat{z}_u^i(k+l|k)$, and $\hat{\sigma}^i_d(k+l+1|k)$ are computed by the following dynamics considering the optimal input trajectory $\delta\hat{\alpha}^i(k:K+N-1|k)$:

\begin{subequations}\label{eq:costrsdynamic}
\begin{align}
&\left\{\begin{array}{lcl}
\hat{\alpha}^i(k+l|k)&=&\hat{\alpha}^i(k+l-1|k)+\delta\hat{\alpha}^{i}(k+l|k)\\
\hat{\alpha}^i(k-1|k)&=&\alpha^i(k-1)
\end{array} \qquad \forall l=0,\cdots,N-1,\right. \label{eq:hat_alpha}
\\
&\left\{\begin{array}{lcl}
\hat{\varepsilon}^i_d(k+l+1|k)&=&\Phi_{ii}\hat{\varepsilon}^i_d(k+l|k)+\Gamma_{ii}\delta\hat{\alpha}^i(k+l|k)\\
\hat{\varepsilon}^i_d(k|k)&=&{\varepsilon}^{i}_{d}(k)
\end{array} \right.\label{eq:hat_eps}
\\
&\left\{\begin{array}{lcl}
\hat{z}_u^i(k+l+1|k)&=&\Phi_{ii}\hat{z}_u^i(k+l|k)+\sum_{j\in\mathcal{N}_{IN}^i}\Phi_{ij}\hat{z}_c^j(k+l|k)+\Gamma_{ii}\left(y_r^i+\hat{\alpha}^i(k+l-1|k)\right)\\
\hat{z}_u^i(k|k)&=&{z}^i(k)-{\varepsilon}^{i}_{d}(k)
\end{array} \right.\label{eq:hat_x_c}
\\
&\left\{\begin{array}{lcl}
\hat{\sigma}^i_d(k+l+1|k)&=&\Phi_{ii}\hat{\sigma}^i_d(k+l|k)+\sum_{j\in\mathcal{N}_{IN}^i}\Phi_{ij}\hat{\sigma}^j_d(k+l|k)+\Gamma_{ii}\Delta\hat{\alpha}^i(k+l)\\
\hat{\sigma}^i_d(k|k)&=&0
\end{array} \right.\label{eq:hat_sigma}
\end{align}
\end{subequations}
\noindent where
\begin{equation}
\label{eqDeltaalpha}
\Delta\hat{\alpha}^{i}(k+l|k)= y_r^{i}(k)+\delta\hat{\alpha}^{i}(k+l|k)-(y_r^{i}(k-1)+\delta\hat{\alpha}^{i}(k+l|k-1))
\end{equation}
denotes the difference between optimal input trajectories at time step $k+l$ predicted at last time step $k-1$ and current time step $k$. Hence, \eqref{eq:hat_sigma} computes the effect of $\Delta\hat{\alpha}^{j}(k+l|k)$ for all $j\leq i$ transmitted to any subsystem $m\in\mathcal{N}_{OUT}^i$ and modifies their optimal state trajectory predicted at time step $k-1$. Therefore, in view of the recursive feasibility of $\delta\hat{\alpha}^{m}(k:k+N-1|k-1)$ corresponding to any subsystem $m\in\mathcal{N}_{OUT}^{i}$ at time-step $k$, it is crucial to predict the value of state $\hat{\sigma}^{m}_{d}(k+l|k)$ using equation \eqref{eq:hat_sigma} for any $l=1,\cdots,N$ and for all $m\in\mathcal{N}_{OUT}^{i}$.

Now, according to the receding horizon principle, the RG control law is set as
\begin{equation}\label{claw}
\alpha^i(k)=\alpha^i(k-1)+\delta\hat\alpha^i(k|k).
\end{equation}

It worth to note that the prediction horizon $N$ and the symmetric positive definite weighting matrices $Q_{ii}$, $P_{ii}$, $P_{\alpha_{ii}}$ and $R_{\alpha_{ii}}$ are free design parameters. In order to obtain the stability results presented below, $P_{ii}$ has to be selected as the positive definite solution of the Lyapunov equation
\begin{equation}\label{eq11}
\left[\Phi_{ii}\right]^TP_{ii}\Phi_{ii}-P_{ii}=-Q_{ii}
\end{equation}
and $P_{\alpha_{ii}}$ can be chosen as a positive definite matrix such that
\begin{equation}\label{eq12}
P_{\alpha_{ii}}>\left[\Gamma_{ii}\right]^TP_{ii}\Gamma_{ii}+R_{\alpha_{ii}}.
\end{equation}
When solving the $i$-RHOP problem with the initial state $z_c^i(k)=z^i(k)$ (the measured state), the optimal solution $\delta\hat{\alpha}^{j}(k|k)$, $j\in \mathcal{N}_{IN}^{i}$ is already available, as well as the future predicted state trajectories $ z_{c}^{j}(k+j|k)$, $j\in \mathcal{N}_{IN}^{i}$, $l=0,...,N-1$, which can be considered as known inputs. The overall proposed algorithm is summarized in Algorithm 2.

\begin{algorithm}\label{AlRHRG}
	\caption{DCT-DRG Algorithm}
	\begin{algorithmic}
		\item[1)] \textbf{Offline computation}
		\begin{itemize}
			\item[1-1)] Compute $\mathbb{X}\mathbb{U}^i(k)$, $k=1,...,N-1$ and $\mathbb{X}\mathbb{U}^i_{\infty}$, for any subsystem $i=1,\cdots,M$ (see Subsection \ref{ss31}).
			\item[1-2)] Compute MOASs $\mathbb{O}^i_{\varepsilon}$ corresponding to any subsystem $i=1,\cdots,M$ (see Subsection \ref{ss32}).
		\end{itemize}
		
		\item[2)] \textbf{Distributed receding horizon online computation in $i$-subsystem at any time step $k$}\\	
		\begin{itemize}	
			\item[2-1)] Measure the current state value $z^i(k)$ of any subsystem and set $z_c^i(k|k)=z^i(k)$. 
			\item[2-2)] Considering $\check{g}^j(k+l|k)$ and $\check{g}^m(k+l|k-1)$ for all $j\in\mathcal{N}_{IN}^{i}$, $m\in\mathcal{N}_{OUT}^{i}$, and $l=0,\cdots,N-1$, solve $i$-RHOP optimization problem in Eq. \eqref{eq18}.
			\item[2-3)] Compute $\alpha^i(k)$ by Eq. \eqref{claw} and $\check{g}^i(k)$ by Eq. \eqref{reference} determining the current admissible reference input.
			\item[2-4)] Return to Step 2-1 at the next time step $k=k+1$.
		\end{itemize}
	\end{algorithmic}
\end{algorithm}
\begin{remark}	\label{Remark:R0}
	At time step $k$, the $i$-RHOP is computing new admissible reference input sequence to apply to the subsystem $i$ and, simultaneously, subsystem $i$ is affecting the subsystems $m\in\mathcal{N}_{OUT}^{i}$ as its outlet neighbors. Therefore, it compromise the feasibility of reference input sequences computed for subsystems $m\in\mathcal{N}_{OUT}^{i}$ at the previous time step $k-1$ and hence, the feasibility of a distributed algorithm is not guaranteed. To overcome this obstacle, as emphasized in \cite{garone2017reference}, the previously presented distributed reference/command governors enforce the updating value of any reference input to belong to a specified static constraint set. However, determining this set is problematic and it, obviously, makes the algorithm to be more conservative. Avoiding such a static limit on the updating value, in this paper dynamic interaction constraint sets \eqref{eq:terminal_sigma} and \eqref{eq:constr_aniticipative} are introduced to preserve the recursive feasibility of the proposed distributed algorithm and allow the $\alpha^i(k)$ to be updated by the possible value without any static limit and in such a way that, the previously computed reference input sequences of the outlet neighbors at time step $k-1$ remain feasible, also, at the current time step $k$.
\end{remark}
\begin{remark}	\label{Remark:R1}
	Comparing to the conventional -one step horizon- reference governors, which is named here as "static reference governor", and in view of \eqref{eq:hat_alpha}, the proposed algorithm enjoys of a "dynamic reference governor strategy". This property, not only increases the degrees of freedom in computation of manipulated reference input leading to a faster transient performance of RG and closed-loop system, but also makes it possible to apply the presented dynamic constraint tightening idea which leads to less conservative distributed procedure. This fact is illustrated using the application example presented in Section \ref{SimulationExample} illustrating the performance of the presented DCT-DRG.
\end{remark}

\begin{remark}	\label{Remark:R2}
As a more conservative but trivial alternative method, one can substitute the transient constraint sets $\mathbb{X}{\mathbb{U}}^{i}(.)$ and $\mathbb{X}{\mathbb{U}}^{m}(.)$ in \eqref{eq:constr_stage} and \eqref{eq:constr_aniticipative} with the steady-state tightened constraint set $\mathbb{X}\mathbb{U}^{i}_{\infty}$ and $\mathbb{X}\mathbb{U}^{m}_{\infty}$, respectively, and also use the evolutionary state of CNM \eqref{CNM} as the initial state of Algorithm 2 in place of the measured state value of any subsystem, i.e. $z_c^i(k|k)=z_c^i(k)$. Henceforth, we call this alternative approach as the Static Constraint Tightening based Distributed Reference Governor (STC-DRG) for comparison with the results of proposed DCT-DRG. In fact, the MPC-based RG aprpoach in which is based the proposed DCT-DRG does not uses dynamic tighthening and thus would lead to the STC-DRG approach.
\end{remark}

\subsection{Properties of the Proposed Approach} \label{ss34}
Now, the recursive feasibility, stability and convergence of the proposed DCT-DRG are stated in the following two theorems.
\begin{theorem}[Recursive feasibility]	\label{T1}
	For any subsystem $i$ considering Assumption \ref{A1}, if at time $k=0$ the optimization problem (\ref{eq18})-\eqref{eq:costrsdynamic} is feasible then it remains feasible $\forall k>0$.
\end{theorem}

\noindent\textbf{Proof:}\\
Assume that the states $ z_u^i(k)$ and $\varepsilon_d^i(k)$ are admissible for the corresponding $i$th-RHOP, whose solution is $\delta\hat\alpha^i(k:K+N-1|k)$.
If the RHOP starts sequentially from the first subsystem to $M$-th subsystem, then it is straightforward to show that in $i$-th-RHOP at time step $k+1$ the $j$-th-RHOP for all $j\in\mathcal{N}_{IN}^{i}$ are feasible. Let define the candidate feasible solution to $i$-th-RHOP as $\delta\hat\alpha^i(k+1:K+N-1|k)$. Then, regarding the fulfilled constraint \eqref{eq:constr_aniticipative} while solving $j$-th-RHOP for all $j<i$ and in view of $\mathbb{X}{\mathbb{U}}^{i}(l+1)\subseteq \mathbb{X}{\mathbb{U}}^{i}(l)$, one can write:
\begin{eqnarray}\label{proof_constr_aniticipative}
H_{ii}(\hat{z}_c^i(k+l|k)+\hat{\sigma}^{i}_{d}(k+l|k+1))\in \mathbb{X}{\mathbb{U}}^{i}(l),&& l=1,...,N-2,
\end{eqnarray}
where, $\hat{\sigma}^{i}_{d}(k+l|k+1)$, which is derived by using \eqref{eq:hat_sigma}, indicates variations of the state sequences of subsystems $j$, $\forall j\in\mathcal{N}_{IN}^{i}$ which are transmitted to subsystem $i$. Moreover, these transmitted variations affects the state of $i$-th subsystem at time step $k+1$ while the input sequence calculated at time step $k$ is frozen.
Thus, fulfillment of the constraint \eqref{eq:terminal_sigma} in $j$th-RHOP $\forall j\in \mathcal{N}_{IN}^{i}$ leads to:
\begin{eqnarray}\label{eq:proof_terminal_sigma1}
\sum_{j\in \mathcal{N}_{IN}^{i}}\Phi_{ij}(\hat{ z}_{c}^{j}(k+N-1|k)+\hat{\sigma}^{j}_{d}(k+N-1|k+1))
\in \mathbb{W}_{ z}^{i},
\end{eqnarray}
which implies that $\hat{z}_c^i(k+N|k+1)\in \mathbb{O}^{i}_{\varepsilon}$. Consequently, the constraints \eqref{eq:constr_stage} are verified for $l=1,\cdots,N$ at time step $k+1$ in view of the feasibility of \eqref{eq18} at step $k$.\\
Considering the candidate input sequence $\delta\hat\alpha^{m}(k+1:k+N-1|k)=(\delta\hat{\alpha}^{m}(k+1|k),\dots, \delta\hat{\alpha}^{m}(k+N-1|k),0)$ for $\forall m\in\mathcal{N}_{OUT}^{i}$ at time step $k+1$ and assume the feasibility of $j$-subsystem $\forall j\in\mathcal{N}_{IN}^{i}$, and also in view of \eqref{eq:hat_sigma}, it is easy to show that:
\begin{eqnarray}\label{constr_aniticipative2}
H_{mm}(\hat{z}_{c}^{m}(k+l|k)+\hat{\sigma}^{m}_{d}(k+l|k+1))\in \mathbb{X}{\mathbb{U}}^{m}(l+1),&& \forall m\in\mathcal{N}_{OUT}^{i}, l=1,\cdots,N-2,
\end{eqnarray}
\begin{eqnarray}\label{proof_terminal_sigma2}
\sum_{j\in\mathcal{N}_{IN}^{m}}\Phi_{mj}(\hat{ z}_{c}^{j}(k+N|k)+\hat{\sigma}^{j}_{d}(k+N|k+1))
\in \mathbb{W}_{ z}^{m},&&  \forall m\in\mathcal{N}_{OUT}^{i},
\end{eqnarray}
which implies that the constraints \eqref{eq:constr_aniticipative} and \eqref{eq:terminal_sigma} are also fulfilled for $i$-th-RHOP at time step $k+1$. Furthermore $\hat{\alpha}^{i}(k+N|k)=\hat{\alpha}^{i}(k+N-1|k).$
Recalling \eqref{eq:exp_out}:
$$c^{i}_{c}(k+N|k)=\mathcal{C}^{i}\begin{bmatrix}{ z}_c^i(k+N|k)\\{y}_r^i+\hat{\alpha}^i(k+N-1)
        \end{bmatrix}$$
and $\eqref{eq:constr_term}$ and the fact that $\mathbb{O}^j_{ z}$ is output-admissible (see \eqref{eq16}), then $c^{i}_{c}(k+N|k)\in\mathbb{X}{\mathbb{U}}_{\infty}^i\subseteq\mathbb{X}{\mathbb{U}}^i(N-1)$. Furthermore, recalling \eqref{eq:exp_state} and $\hat{z}_c^i=\hat{z}_{u}^{i}+\hat{\varepsilon}^{i}_{d}$ yields to:
$$\begin{array}{ll}
\begin{bmatrix}{z}_c^i(k+N+1|k)\\{y}_r^i+\hat{\alpha}^i(k+N)
        \end{bmatrix}&=\begin{bmatrix}{ z}_c^i(k+N+1|k)\\{y}_r^i+\hat{\alpha}^i(k+N-1)
        \end{bmatrix}\\
        &=\mathcal{A}^{i}\begin{bmatrix}{ z}_c^i(k+N|k)\\{y}_r^i+\hat{\alpha}^i(k+N-1)
        \end{bmatrix}+\mathcal{B}^iw_{ z}^i(k+N|k),\end{array}$$
since \eqref{eq:constr_term} holds at $k$ for all $j\leq i$, then $\hat{ z_{c}}^{j}(k+N|k)\in \mathbb{O}^j_{ z}$ is fulfilled for all $j\in \mathcal{N}_{IN}^{i}$. Therefore $w_{ z}^i(k+N|k)=\sum_{j\in \mathcal{N}_{IN}^{i}}\Phi_{ij}\hat{ z_{c}}^{j}(k+N|k)\in \mathbb{W}^i_{ z}$ and, in view of the fact that $\mathbb{O}_{\varepsilon}^i$ is RPI, then $({ z}_c^i(k+N+1|k),{y}_r^i+\hat{\alpha}^i(k+N))\in \mathbb{O}_{\varepsilon}^i$, and therefore constraint \eqref{eq:constr_term} is also verified at instant time step $k+1$.\hfill$\square$

\begin{theorem}[Stability and convergence] \label{T2}
    Assume that $\varepsilon_d^{j}(k)$ $\forall j\in\mathcal{N}_{IN}^{i}$ has converged to zero, the optimization problem (\ref{eq18})-\eqref{eq:hat_eps} is feasible for all time steps, and the Assumption \ref{A1} holds. Then, the control law \eqref{claw} steers the output $y^{i}$ to the admissible set point $y^i_r+ \alpha^i_{ad}$, where
        \begin{align}
        \label{eq:aad_def}
        \alpha^i_{ad}=\mathop{\text{argmin}}_{H_{ii}\left(I-\Phi_{ii}^{-1}\Gamma_{ii} (y^{i}_r+\alpha^i)\right) \in \mathbb{X}_{\mathbb{U}}^{i}(\varepsilon)}\|\alpha^i\|^2_{P_{\alpha_{ii}}}
        \end{align}
    respects the state and input constraints robustly.
\end{theorem}

\noindent\textbf{Proof:}\\
\textit{i) Stability and convergence to zero of $\varepsilon_d^{i}(k)$}:\\
While the convergence of $i$-th subsystem is being investigated, the assumption of convergence of $\varepsilon_d^{j}(k)$ for all $j\in \mathcal{N}_{IN}^{i}$ to zero make sense due to inductive reasoning for all priori subsystems starting from $j=1$.
Consider the feasible and possibly suboptimal solution of \eqref{eq18} for subsystem $i$ at time step $k+1$, i.e., $\delta\hat\alpha^{i}(k+1:K+N-1|k)$. Sub-optimality means:
$$J_{N}^{i^{*}}(\varepsilon_{d}^{i}(k+1),y^i_r,\alpha^i(k))\leq
J_{N}^{i}(\delta\hat\alpha^{i}(k+1:K+N-1|k);
\varepsilon_{d}^{i}(k+1),y^i_r,\alpha^i(k)).
$$
Considering convergence of all subsystems $j\in\mathcal{N}_{IN}^{i}$, for the subsystem $i$ and according to standard arguments in RHC control, it follows that
$$\begin{array}{ll}J^{i}_{N}\left(\delta\hat\alpha^{i}(k+1:k+N-1|k);
\varepsilon^{i}_{d}(k+1), z_{u}^{i}(k+1),y^i_r,\alpha^i(k)\right)= \\ \qquad
J_{N}^{i^{*}}\left(\varepsilon^{i}_{d}(k),y^i_r,\alpha^i(k-1)\right)
-\|\hat{\varepsilon}^{i}_d(k|k)\|^2_{Q_{ii}}-\|\delta\hat{\alpha}^{i}(k|k)\|^2_{R_{\alpha_{ii}}}
\end{array}$$
such that
$$\begin{array}{c}J_{N}^{i^{*}}(\varepsilon^{i}_{d}(k+1),y^i_r,\alpha^i(k))-
J_{N}^{i^{*}}\left(\varepsilon^{i}_{d}(k),y^i_r,\alpha^i(k-1)\right)\leq
-\|{\varepsilon}^{i}_d(k)\|^2_{Q_{ii}}-\|\delta\hat{\alpha}^{i}(k|k)\|^2_{R_{\alpha_{ii}}}.
\end{array}$$
Due to the definite positiveness and non-increasing evolution of the optimal cost function $J_{N}^{i^{*}}$, it yields
\begin{align}
&J_{N}^{i^{*}}\left(\varepsilon^{i}_{d}(k),y^{i}_r,\alpha^i(k-1)\right)\rightarrow \bar{J}^{i}\text{   as   }k\rightarrow+\infty\label{limJ}\\
&\lim_{k\to\infty}\|{\varepsilon}^{i}_d(k)\|^2_{Q_{ii}} =\lim_{k \to \infty}\|\delta\hat{\alpha}^{i}(k|k)\|^2_{R_{\alpha_{ii}}}=0 \label{eq21}
\end{align}
In view of~\eqref{eq18},~\eqref{eq:costrs},~\eqref{limJ}, and~\eqref{eq21}, one can write
\begin{align}
\|\hat{\alpha}^{i}_{k+N-1|k}\|^2_{P_{\alpha_{ii}}}\rightarrow \bar{J}^{i}\text{   as   }k\rightarrow+\infty
\label{eq:limalpha}
\end{align}
Eventually, \eqref{limJ} and \eqref{eq21} imply that ${\alpha}^{i}_{k}\rightarrow \bar{\alpha}^{i}$ as $k\rightarrow+\infty$, where $\|\bar{\alpha}^{i}\|^2_{P_{\alpha_{ii}}}=\bar{J}^{i}$ and $H_{ii}((I-\Phi_{ii})^{-1}\Gamma_{ii}(y^{i}_r+\alpha^i))\in\mathbb{X}_{\mathbb{U}}^{i}(\varepsilon)$, i.e., it is admissible in view of ~\eqref{eq:constr_term}. In steady-state conditions $J^{i^{*}}_N=\bar{J}^{i}$, where the corresponding solution of \eqref{eq18} is given by $\delta\hat\alpha^{i^{ss}}=\begin{bmatrix} 0 & \dots & 0\end{bmatrix}^T$. So, the constant trajectories of $\hat{\varepsilon}^{i}_{d}$ and $\hat{ z}_{u}^{i}$ will be equal to $\hat{\varepsilon}^{i^{ss}}_{d}$ and $\hat{z}_{u}^{i^{ss}}$, respectively.

By contradiction, assume that $\bar{\alpha}^{i}\neq \alpha^i_{ad}$, where $\alpha^i_{ad}$ is defined in~\eqref{eq:aad_def}. Then, an alternative solution to the problem \eqref{eq18} is $\delta\tilde\alpha^{i}=\begin{bmatrix} 0 & \dots & \delta\alpha^i\end{bmatrix}^T$ where $\delta\tilde{\alpha}^{i}=\lambda(\alpha^i_{ad}-\bar{\alpha}^{i})$ and $\lambda\in(0,1)$.\\
By convexity of $\mathbb{X}^{i}_\mathbb{U}(\varepsilon)$, $\tilde{\alpha}^{i}=\bar{\alpha}^{i}+\delta\tilde{\alpha}^{i}$ is admissible. The corresponding trajectory of $\hat{ z}_{u}^{i}$ is constant and, during the whole prediction horizon, is equal to $\hat{ z}_{u}^{i^{ss}}$. On the other hand, $\hat{\varepsilon}^{i}_d(k+l)=0$ for $l=0,\cdots,N-1$ and $\hat{\varepsilon}^{l}_d(k+N)=\Gamma_{ii} \delta\tilde{\alpha}^{i}$, which is feasible for a sufficiently small value of $\lambda$, since $\varepsilon>0$. So, the value of $J^{i}_N$ which is computed according to \eqref{eq18} with this alternative solution is equal to:
\begin{equation}
\begin{array}{lcl}
\tilde{J}^{i}_N&=&\|\lambda\Gamma_{ii} (\alpha^i_{ad}-\bar{\alpha}^{i})\|^2_{P_{ii}}+\|\lambda(\alpha^i_{ad}-\bar{\alpha}^{i})\|^2_{R_{\alpha_{ii}}}+
\|(1-\lambda)\bar{\alpha}^{i}+\lambda \alpha^i_{ad}\|^2_{P_{\alpha_{ii}}}\\
&=&\lambda^2\|\bar{\alpha}^{i}-\alpha^i_{ad}\|_{\Gamma_{ii}^TP_{ii}\Gamma_{ii}+R_{\alpha_{ii}}}+
(1-\lambda)^2\|\bar{\alpha}^{i}-\alpha^i_{ad}\|^2_{P_{\alpha_{ii}}}+\|\alpha^i_{ad}\|^2_{P_{\alpha_{ii}}}\\
&&+2(1-\lambda)\alpha^{i^T}_{ad} P_{\alpha_{ii}} (\bar{\alpha}^{i}-\alpha^i_{ad}).
\end{array}
\label{eq:Jalt}
\end{equation}
Note also that
$$\bar{J}^{i}=\|\bar{\alpha}^{i}\|^2_{P_{\alpha_{ii}}}=\|\bar{\alpha}^{i}-\alpha_{ad}^{i}\|^2_{P_{\alpha_{ii}}}+
\|\alpha^i_{ad}\|^2_{P_{\alpha_{ii}}}+2\alpha^{i^T}_{ad} P_{\alpha_{ii}} (\bar{\alpha}^{i}-\alpha^i_{ad})$$
Therefore, since $P_{\alpha_{ii}}$ verifies \eqref{eq12}, it could be shown that
$$\bar{J}^{i}-\tilde{J}^{i}_N\geq \|\bar{\alpha}^{i}-\alpha^i_{ad}\|^2_{P_{\alpha_{ii}}}(1-\lambda^2+(1-\lambda)^2)+2\lambda \alpha^{i^T}_{ad} P_{\alpha_{ii}} (\bar{\alpha}^{i}-\alpha^i_{ad})$$
we have $(1-\lambda^2+(1-\lambda)^2)>0 \quad \forall\lambda\in(0,1)$ and $2\alpha_{ad}^{i^{T}} P_{\alpha_{ii}} (\bar{\alpha}^{i}-\alpha^i_{ad})\geq 0$ due to optimality of $\alpha^i_{ad}$ with respect to the quadratic function $\|\alpha^i\|^2_{P_{\alpha_{ii}}}$ in the admissible set. Therefore
$$\bar{J}^{i}>\tilde{J}^{i}_N$$
which contradicts the assumption that $\bar{\alpha}^{i}\neq \alpha^i_{ad}$ corresponds to a steady-state for the subsystem $i$ controlled with the $i$-th-RHOP control law. Therefore, the only steady-state, compatible with \eqref{eq18}, is the one corresponding to the condition $\bar{\alpha}^{i}= \alpha^i_{ad}$.\\
\\
\textit{ii) Convergence to the reference}:\\
Consider the whole large scale system \eqref{eq1} at the steady state. Then, according to the closed-loop nominal model \eqref{CNM} and the fact that $ z_c^i(k+1)= z_c^i(k)= z_{c_{ss}}^{i}$ for $i=1,\cdots,M$, leads to
\begin{equation}
\label{SSCCNM}\begin{array}{lcl}
(I-\Phi_{ii}) z^{i}_{c_{ss}}&=& \sum_{j\in \mathcal{N}_{IN}^{i}}\Phi_{ij} z^{j}_{c_{ss}}
+\Gamma_{ii}\tilde{y}_r^i.  \\
\end{array}
 \end{equation}
and following the conditions in Assumption \ref{A1}, it yields:
$$
y_{c_{ss}}^{i}=\Upsilon_{ii} z^{i}_{c_{ss}}=\tilde{y}_r^i
$$
The same reasoning can be applied to all the subsystems closed-loop sub-system \eqref{InterconnectedSubsystem}. Thus, in the case of immeasurable asymptotically constant and bounded disturbance signals entering to each sub-system, if follows $y^{i}_{ss}=\Upsilon_{ii} z^{i}_{ss}=\tilde{y}_r^i$. And, in the case of not measurable varying and bounded disturbance sequences which affect each sub-system and applying the first sample of the manipulated reference sequence in any time step, we have $y^{i}(k\rightarrow\infty)\in\Upsilon_{ii}( z^{i}_{c_{ss}}\oplus\Omega_{ii}\mathbb{W}_{e}^{i})=\tilde{y}_r^i\oplus\Upsilon_{ii}\Omega_{ii}\mathbb{W}_{e}^{i}$\hfill$\square$\\

\section{Case study} \label{SimulationExample}

\subsection{Description}

\begin{figure}[!h]
	\begin{center}
		\includegraphics[width=15cm]{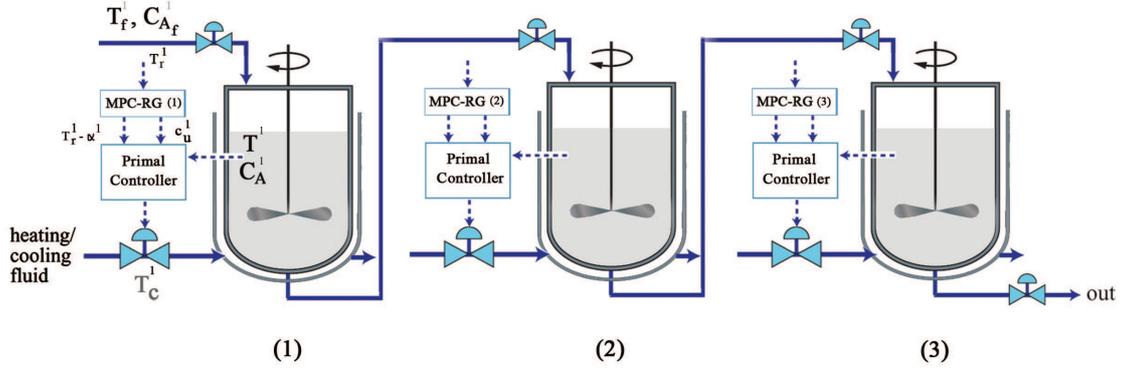}
		\caption{Schematic diagram of three serially connected CSTR}
		\label{fig:aa}
	\end{center}
\end{figure}

Consider the system in Figure \ref{fig:aa} that is created as cascade connection of three jacketed Continuous Stirred Tank Reactors (CSTR).  The model of this systme is obtained considering that a single irreversible reaction $A \rightarrow B$ in each CSTR is represented by  \cite{mclain2000nonlinear,magni2001stabilizing,ferramosca2009mpc}:
$$
\begin{array}{lcl}
\dot{C}_{A} &=& \frac{q}{V}(C_{Af}-C_{A})-k_{0}exp(-\dfrac{E}{RT})C_{A} \\
\dot{T} &=& \frac{q}{V}(T_f-T)+\frac{\Delta H}{\rho C_{p}}k_0 exp(-\dfrac{E}{RT})C_{A}+\dfrac{UA}{V\rho C_{p}}(T_c-T) \\
y &=& \begin{bmatrix}
0 & 1
\end{bmatrix} \begin{bmatrix}
C_A \\
T
\end{bmatrix}
\end{array}
 $$
where, $C_{A}$ and  $T$ are reactant concentration and reactor temperature, respectively. Also, the cooling liquid temperature $T_{c}$ is the input variable. Nominal values of other parameters and coefficients are reported in Table \ref{Table:T1}.

\begin{center}
\begin{table}[h]
\caption{Nominal values of CSTR variables}
{\begin{center}
\begin{tabular}{|cc|cc|}
\hline
Variable & Nominal value & Variable & Nominal value \\
\hline
$C_{p}$ & $0.239 J/g.K$ & $\Delta H$ & $5\times 10^{4} J/mol$ \\
$q$ & $100 L/min$ & $k_{0}$ & $7.2 \times 10^{10} min^{-1}$ \\
$V$ & $100 L$ & $UA$ & $5\times 10^{4} J/min.K$ \\
$E/R$ & $8750 K$ & $\rho$ & $1000 g/L$ \\
\hline
\end{tabular}
\end{center}}
\label{Table:T1}
\end{table}
\end{center}

This model is linearized around the operating point $C_{Af}=1 mol/l$, $T_{f}=300 K$, $T_{c}=300 K$, $\bar{T}=301.15 K$, $\bar{C_{A}}=0.98296 mol/l$ and discretized with sampling time $\Delta T_s=0.6 min$. Consequently, for each subsystem $i$, linear discrete time model is described as:
$$
\begin{array}{lcl}
\begin{bmatrix}
\Delta C_{A}^{i}(k+1) \\
\Delta T^{i}(k+1)
\end{bmatrix}&=&\begin{bmatrix}
0.54271 & -3e^{-4} \\
0.73488 & 0.19196
\end{bmatrix}\begin{bmatrix}
\Delta C_{A}^{i}(k) \\
\Delta T^{i}(k)
\end{bmatrix}+\begin{bmatrix}
-3e^{-4} \\
0.6152
\end{bmatrix}\Delta T_c^{i}(k)+\begin{bmatrix}
0.2 & 0 \\
0 & 0.2
\end{bmatrix}\begin{bmatrix}
\Delta C_{A}^{j}(k) \\
\Delta T^{j}(k)
\end{bmatrix}+w^i(k) \\
\Delta y(k)&=&\begin{bmatrix}
0 & 1
\end{bmatrix}\begin{bmatrix}
\Delta C_{A}^{i}(k) \\
\Delta T^{i}(k)
\end{bmatrix},\ \  i=1,2,3,\ \   \mathcal{N}_{IN}^{i>1}=i-1,
\end{array}
$$
where:
\begin{itemize}
\item $\vert w^i\vert \leq [0.05, 0.5]^T, i=1,2,3$ are unknown and but bounded state disturbances and 
\item $[\Delta C_{A}^{i}(0), \Delta T^{i}(0)]^T =[0,0]^T$.
\end{itemize}


The control goal is that the reactor temperature $\Delta T^{i}$ track the set-point $\Delta T_r^{i}$ or the nearest feasible value (in case of the inadmissible $\Delta T_r^{i}$ regarding the system constraints) without any offset in the presence of disturbance. The constraints to be robustly satisfied at each time instant are $-3\leq\Delta T_c^i\leq 3,-5\leq\Delta T^{i}\leq 5$, $i=1,2,3$.
\begin{center}
\begin{table}[h]
\centering
\caption{Disturbance sequences entering the subsystems}\resizebox*{9cm}{!}
{\begin{tabular}{|c|ccc|}
\hline
 & $8<k\leq 100$ & $100<k\leq 125$ & $k>125$\\
\hline
 $w^{i}, i=1,2,3$ & $\begin{bmatrix}
-0.05 \\
0.5
\end{bmatrix}$ & $ \begin{bmatrix}
0.05 \\
-0.5
\end{bmatrix}$& $\begin{bmatrix}
0.05 \\
0.5
\end{bmatrix} rand(1,1)$ \\
\hline
\end{tabular}}
\label{Table:T3}
\end{table}
\end{center}

In order to satisfy Assumption \ref{A1}, following Lemma 2.1, local decentralized tracking LQR controllers are implemented \cite{shahram2013mpc}. Each linear model was augmented by an error integrator dynamic and then according to the approach of LQR, stabilizing state feedback are computed for each locally augmented model.  Transient and steady state constraints and, also, MOASs are computed using the approaches and algorithms presented in Sections \ref{ss31} and \ref{ss32}, respectively. 

For evaluating the performance of proposed approach, two controllers are considered: 
\begin{itemize}
\item the first controller follows proposed DCT-DRG approach, while
\item the second controller is following the SCT-DRG approach recalled in Remark \ref{Remark:R2}. 
\end{itemize}

For both controllers and for any RHOP \eqref{eq18}, the matrices $Q_{ii} = I$, $R_{ii}=1$, and $R_{\alpha_{ii}}=1$ are considered and the corresponding matrices $P_{ii}$ and  $P_{\alpha_{ii}}=2(\Gamma'P^{ii}\Gamma+R_{\alpha})$, are locally computed for $i=1,2,3$. The simulated results are obtained by the prediction horizon $N=3$. 

\subsection{Results}

\begin{figure}[H]
	\begin{center}
		\includegraphics[width=17cm]{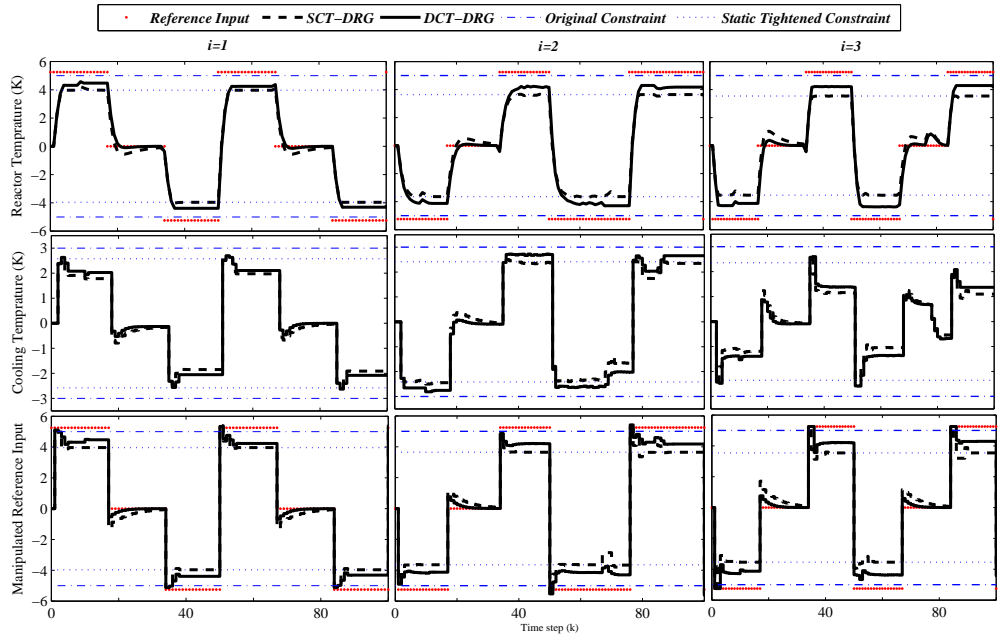}
		\caption{Closed-loop response, Input signal, and manipulated reference input corresponding to each subsystem for SCT-DRG ( \textit{dashed line}) and DCT-DRG ( \textit{solid line})}
		\label{fig:sys1}
	\end{center}
\end{figure}

Figure \ref{fig:sys1} shows simulation results for subsystems 1 to 3. Each column is corresponding to one subsystem which has been determined by subsystem index. The three rows from top to bottom, illustrates the output, input, and manipulated reference input, respectively, provided by the proposed DCT-DRG (solid line) and SCT-DRG (dashed line). The original and steady-state tightened constraint borders regarding the input and state (output) constraints have been also depicted with dashed-dotted and dashed blue lines, respectively. The red dotted line indicate the original reference input. 

\begin{figure}[H]
	\begin{center}
		\includegraphics[width=16cm]{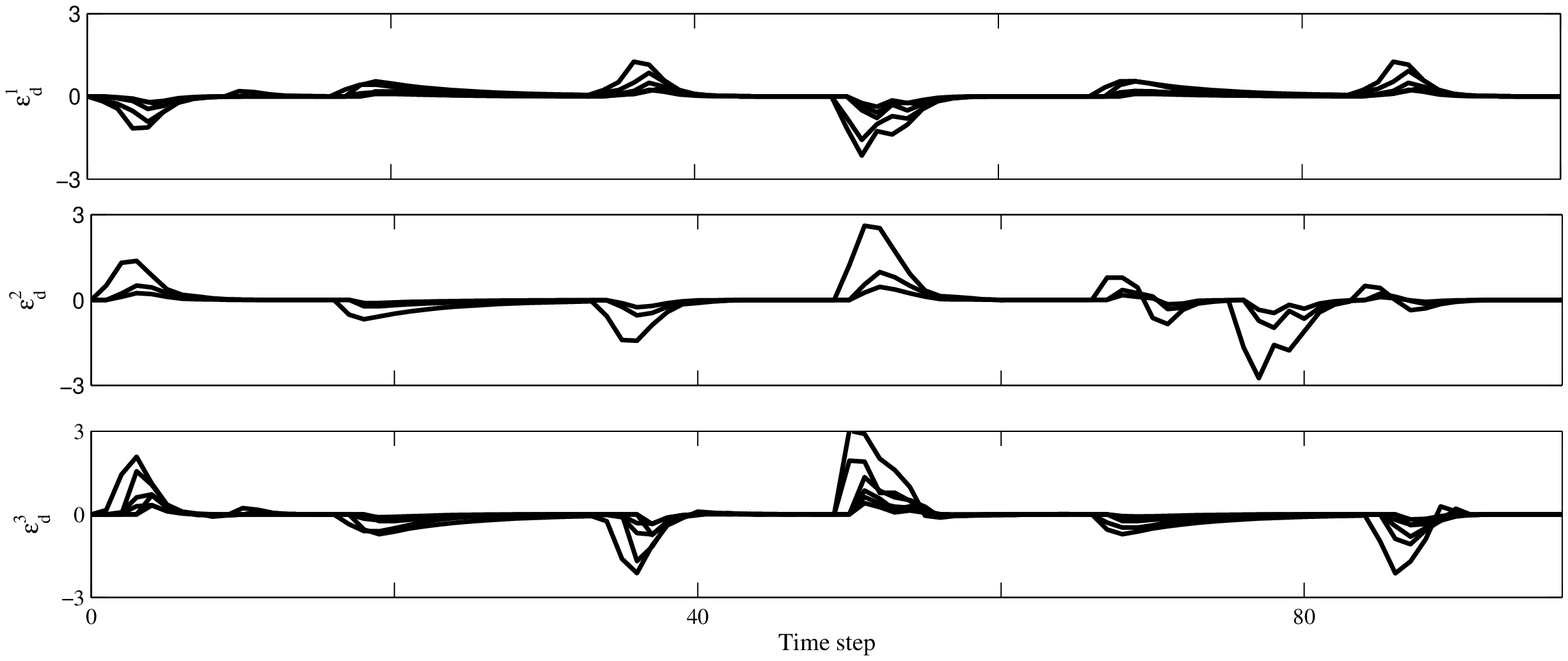}
		\caption{controlled error variable $\epsilon_d^i$ for $i=1,2,3,$ according to CEM (\ref{EM})}
		\label{fig:sys2}
	\end{center}
\end{figure}

The top and middle diagrams in Figure \ref{fig:sys1} presents the output and control input diagrams, respectively. The more desirable, less conservative, and reliable closed-loop performance of the proposed DCT-DRG is achieved in comparing the scenario of SCT-DRG. The results also show, clearly, a better performance in converging to admissible set-points when using DCT-DRG. The robust constraint satisfaction is achieved also for inadmissible references. The bottom charts of Figure \ref{fig:sys1} presents the manipulated reference inputs $\Delta T_r^{i}+\alpha_{ad}^i$ for $i=1,2,3$ for both scenarios.\\
Figure \ref{fig:sys2} shows that controlled error variable $\epsilon_d^i$ for $i=1,2,3,$ according to CEM (\ref{EM}), always tend to zero asymptotically and confirms the convergence results of proposed algorithm.

\begin{figure}[H]
	\begin{center}
		\includegraphics[width=14cm]{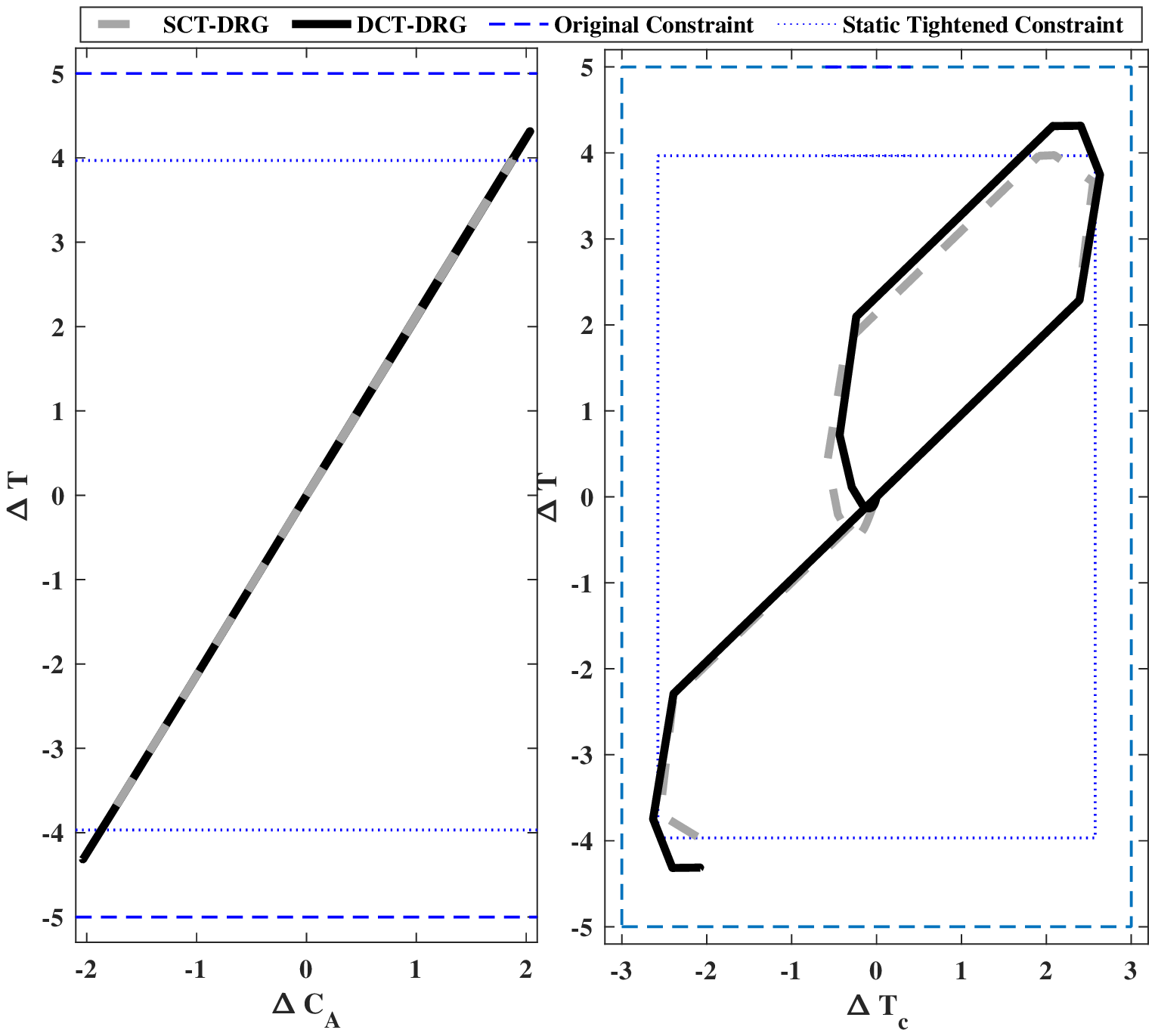}
		\caption{State space response of the first subsystem without disturbance}
		\label{fig:nSS}
	\end{center}
\end{figure}

\begin{figure}[H]
	\begin{center}
		\includegraphics[width=14cm]{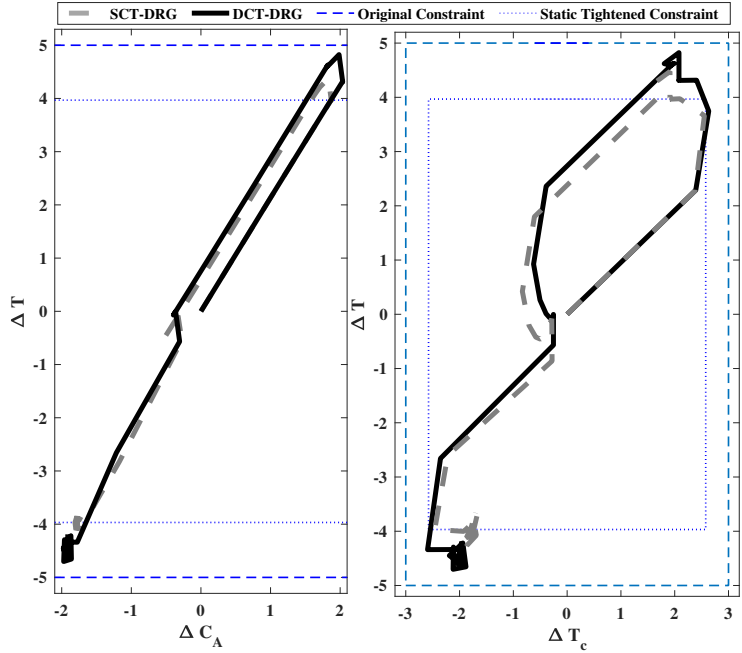}
		\caption{State space response of the first subsystem with disturbance}
		\label{fig:dSS}
	\end{center}
\end{figure}

Figure \ref{fig:nSS} presents for the first subsystem the evolution of the states in the state space when both  DCT-DRG (solid line) and SCT-DRG (dash-line) approaches are used in case no disturbances are applied to system.  Figure \ref{fig:dSS} presents the same evolution in case disturbances are affecting the system. From these figures, it can be seen that the evolution of the states are more constrained in case of SCT-DRG because the original constraints are tightened with the steady steady invariant set. These more tightened constrained justifies the worst results obtained with SCT-DRG compared to those achieved with the DCT-DRG.


\subsection{Discussion}

Regarding the simulation results, the following noticeable facts can be highlighted:
\begin{itemize}
    \item In the case of inadmissible references, in transient and/or steady state, the nearest feasible value is computed by the proposed DRG.
    \item In case of transient response, the proposed algorithm provide a sequence of reference values instead of a single value obtained by the classical RG approach. This property is one of the main advantages of the proposed MPC-RG over the conventional RGs and is more visible with the proposed DCT-DRG. This property of the proposed DRG is preserved due to dynamic interaction constraint set and variable reference prediction horizon as discussed in Remarks \ref{Remark:R0} and \ref{Remark:R1}.
    
    \item The manipulated reference input also imply that, in case of feasible reference inputs, the value of $\delta\alpha^i$ for $i=1,2,3$ tends to zero asymptotically.
    \item The obtained results also show that by applying the SCT-DRG, the output, input, and the steady-state value of manipulated reference input of subsystems are completely limited to the steady-state constraint constraint sets. On the other hand, DCT-DRG leads to a less conservative performance which allows the local RG to determine their manipulated reference input out of their steady-state tightened constraint sets, if needed. Consequently, the local RG are able to manage their input and state (output) out of the steady-state tightened constraint to achieve better performance in case of inadmissible set-points while guaranteeing the robust constraint satisfaction. Such a performance for RG is achieved because of the receding horizon implementation of the proposed DCT-DRG while, the performance of other conventional RG algorithms (\cite{garone2011sensorless}, \cite{garone2017reference} and references therein) will be limited in tightened constraint is the same way as SCT-DRG.

\end{itemize}

%



Consequently, all simulation results show a reliable performance of the proposed distributed algorithm and confirm the results in Theorems \ref{T1} and \ref{T2}, while all of the off-line and on-line computations are done locally.
%
%

\section{Conclusions}\label{conclusion}
In this paper, a hierarchical Dynamic Constraint Tightening Decentralized Reference Governor (DCT-DRG) has been developed for constrained large scale systems described by model that can be decomposes in a Lower Block Triangular (LBT) form, as e.g., the case of cascade systems. Then recursive feasibility, stability, convergence, and robust constraint satisfaction of the proposed algorithm proved. The proposed DCT-DRG performs in hierarchical manner enjoying the receding horizon properties since uses as starting point the MPC-based approach introduced in \cite{shahram2013mpc}. All the invariant and output admissible sets are computed locally off-line. By imposing a dynamic interaction constraint, the recursive feasibility preserved. The algorithm is applied to a system composed of a cascade connection of three jacketed Continuous Stirred Tank Reactors (CSTR). The results confirm the properties such as feasibility, stability, convergence, and robust constraint satisfaction.  
\\\\
\section*{Acknowledgement}
The authors would like to thank Prof. R. Scattolini and Prof. M. Farina for fruitful discussions on the topic.

\bibliographystyle{plain}
\bibliography{Tube_DRG_v31}
\end{document}